\documentclass[10pt]{article}
 
 \usepackage[final]{epsfig}
\usepackage{graphics}
\usepackage{amsmath}
\usepackage{amsfonts}
\usepackage{latexsym}
\usepackage{amssymb}
\usepackage{graphicx}
\usepackage{url}
\newtheorem{lemma}{Lemma}[section]
\newtheorem{proposition}[lemma]{Proposition}
\newtheorem{remark}[lemma]{Remark}

\newtheorem{theorem}{Theorem}

\newtheorem{corollary}[lemma]{Corollary}

\newtheorem{conjecture}[theorem]{Conjecture}

\newcommand{\eps}{{\varepsilon}}

\newcommand{\proofend}{$\Box$\bigskip}

\newcommand{\R}{{\mathbb R}}

\newcommand{\RP}{{\mathbb {RP}}}

\def\proof{\paragraph{Proof.}}

\renewcommand{\v}[1]{\mathbf#1}                         % for boldface vector notation
\newcommand\origin{{\mathcal O}}

\renewcommand\l{\ell}
\renewcommand\a{\alpha}
\renewcommand\b{\beta}

\renewcommand\th{\theta}

\begin{document}

\title{Tractrices, Bicycle Tire Tracks, Hatchet Planimeters, and a 100-year-old Conjecture}

\author{R. L. Foote\footnote{Department of Mathematics and Computer Science, Wabash College, 
Crawfordsville, IN 47933
}  \and   M. Levi\footnote{Mathematics Department, Penn State University, University Park, PA 16802}    \and S. Tabachnikov\footnote{Mathematics Department, Penn State University, University Park, PA 16802}  }
 
\date{}
\maketitle

\section{Introduction} \label{intro}

The geometry of the tracks left by a bicycle has received much attention recently \cite{Dun,Finn1,Finn2,Howe,Levi,Tab}. In this paper we discuss the connection between the motion of a bicycle and that of a curious device known as a hatchet planimeter, and we will prove a conjecture about this planimeter that was made in 1906.

{\bf Bicycle}. We use a very simple model of a bicycle as a moving segment in the plane. The segment has fixed length $\l$, the wheelbase of the bicycle. We denote the endpoints of the segment by $F$ and $R$ for the front and rear wheels. The motion is constrained so that the segment is always tangent to the path of the rear wheel. We will refer to this as the ``bicycle constraint''. This non-holonomic constraint is due to the fact that the rear wheel is fixed on the frame, whereas the front wheel can steer. 
The configuration space of a segment of fixed length is 3-dimensional, and the bicycle constraint defines a completely non-integrable 2-dimensional distribution on it. This is an example of a contact structure, see, e.g., \cite{Arn, Gei}; we shall not dwell on this connection with contact geometry.  

\begin{figure}[hbtp]
\centering
\includegraphics[width=3in]{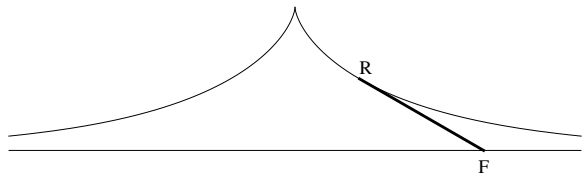}
\caption{The classical tractrix}
\label{tractrix}
\end{figure}

If the path of the front wheel $F$ is prescribed then the rear wheel $R$ follows a constant-distance pursuit curve. The trajectory of the rear wheel is uniquely determined once the initial position of the bicycle is chosen. For example, when $F$ follows a straight line, $R$ describes the classical tractrix, see Figure \ref{tractrix}. More generally, one may call the trajectory of  the rear wheel $R$ the tractrix of the trajectory of the front wheel $F$.

On the other hand, if the path of the rear wheel $R$ is given then the trajectory of the front wheel $F$ is uniquely defined once one fixes the direction of the vector ${RF}$, for which there are two choices. The choice is determined by a coorientation of the rear wheel trajectory, that is, by a continuous choice of a normal direction to it. Given a coorientation, the vector ${RF}$ is determined by the rule that it makes a positive frame with the coorientation. 

One might think that, equally well, one could choose an orientation of the rear wheel track. However, this is not the case: unlike coorientation, orientation is discontinuous at cusps. Cusps are common for the rear-wheel trajectory: they occur when the steering angle equals $90^\circ$, they are where the rear wheel changes its rolling direction. Unless stated otherwise, we assume the trajectory of the front wheel $F$ is smooth.

The study of tractrices goes back to I. Newton (1676), followed by Huygens, Leibniz and Euler. To quote from \cite{Cady}, ``...Euler treated the problem so completely that little or nothing on the subject has appeared since." We hope to make a contribution to this classical subject here.

{\bf Prytz Planimeter}. 
The 19$^\text{th}$ Century was a golden age for mechanical innovation. One modest but very useful device was the planimeter, first invented in Bavaria in 1814.\footnote{The forerunner of the modern bicycle  was invented at about the same time, in 1817, by  Baron Karl von Drais; the invention was  called Draisienne or Laufmaschine.} A planimeter is an instrument that is used to measure the area of a plane figure by tracing around its boundary. As such, it is a mechanical manifestation of  Green's Theorem. There are many types of planimeters, and many improvements have been made over the years, including contributions of Lord Kelvin and James Maxwell. One of the most popular ones  was the polar planimeter, introduced in 1854 by Jacob Amsler, a Swiss mathematician and inventor. See \cite{Bry,Care,Cra,Foote1,Foote2,Hen,Hill,Mur,Ped} for a sampler of the vast literature on planimeters.

By comparison, Amsler's planimeter was more accurate, more compact, and easier to use
than the earlier instruments, and the older ones quickly became obsolete (\cite{Hen} p. 508). Nevertheless, to be accurate the planimeter had to be carefully designed and precisely manufactured, and it could be unaffordable for an engineer of modest means. In the late 1800s, Holger Prytz, a Danish cavalry officer and mathematician, devised an economical and simple alternative to Amsler's planimeter  \cite{Prytz,Ped}.  

\begin{figure}[hbtp]
\centering
\includegraphics[width=1.7in]{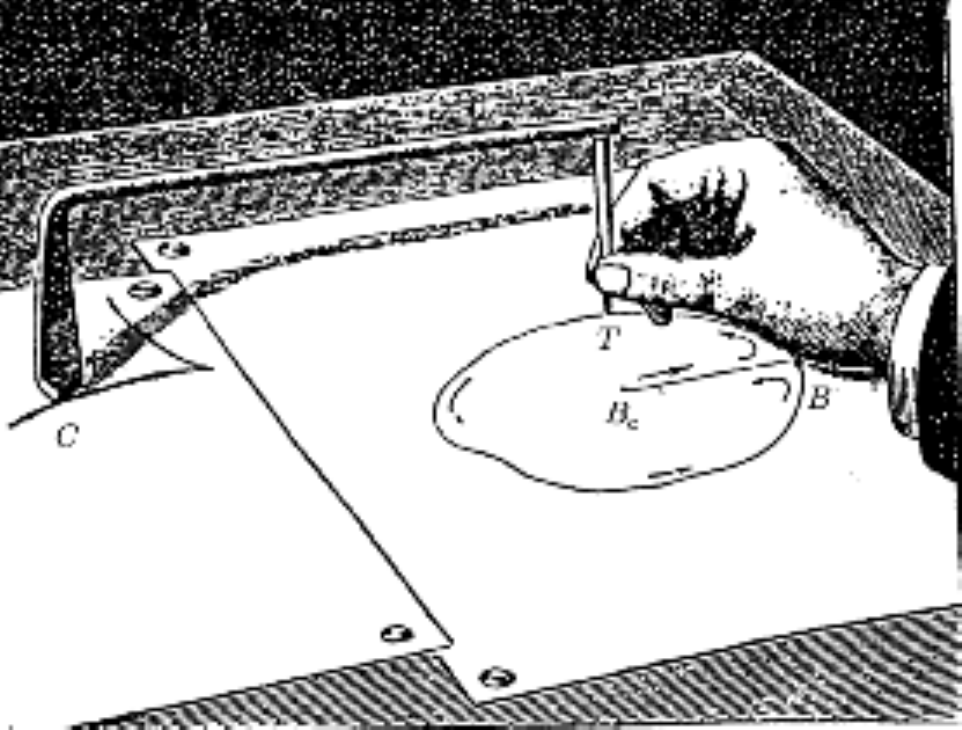}\qquad 
\includegraphics[width=2.3in]{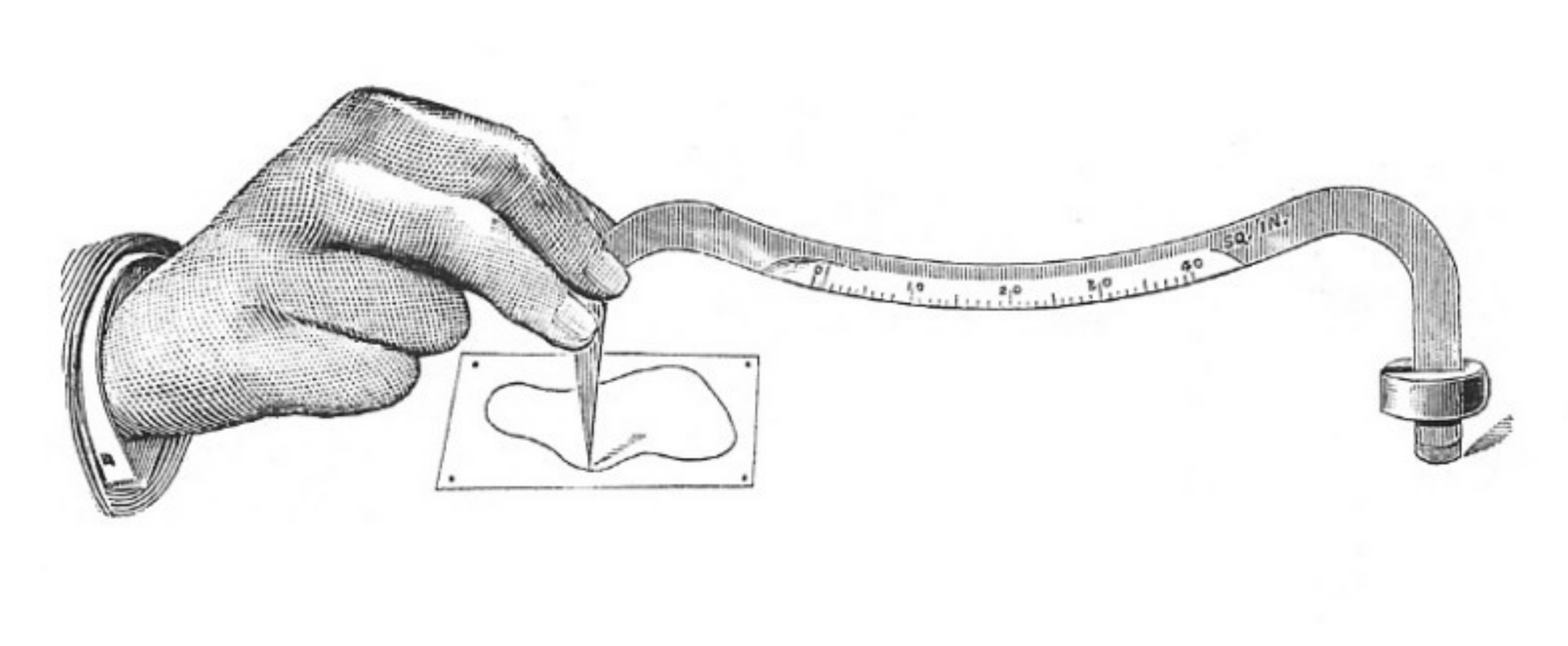}
\caption{The Prytz Planimeter and its modification, due to Goodman}
\label{Prytz}
\end{figure}

Prytz's planimeter consists of a metal rod whose one end, the tracer point, is sharpened to a point, the other end is sharpened to a chisel edge parallel to the rod. (The chisel edge is usually  rounded, making it look similar to a hatchet, and consequently the device is also known as a ``hatchet planimeter'').  It is used  by  guiding the tracer point along a curve, taking care not to impart any torque. The chisel edge tracks along a curve always tangent to the rod. Thus  the hatchet planimeter satisfies the bicycle constraint, with the chisel edge and tracer point playing the roles of the rear and front wheels, respectively. 

\begin{figure}[hbtp]
\centering
\includegraphics[width=2.5in]{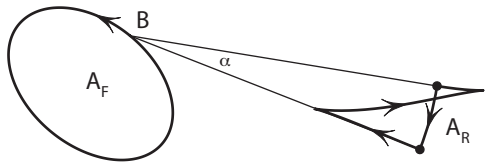}
\caption{Measuring $ A_F $ via   (\ref{diff}).}
\label{planimeter}
\end{figure}

It seems unlikely that something as simple (some would say crude) as a hatchet planimeter could measure area. To use it, put the tracer point at some  point  on the boundary of the region  and  trace around its boundary. The chisel edge follows a zig-zag path, that is, a path with cusps, similar to the rear wheels of a car when
parallel parking, see Figure \ref{planimeter}. If the region is small relative to $\l$, the angular deflection of the planimeter is small.  When the tracer point returns to the initial position, the chisel edge comes to rest in a slightly different position and makes an angle $\alpha$ with its initial position. The area of the region is $\alpha \l^2$, at least approximately. There is an inherent error, which actually makes the hatchet planimeter more mathematically interesting than its exact cousins. As we will see, understanding the source of this error can help the user minimize it. 

{\bf Menzin's conjecture}. The trajectory of the front wheel of a bicycle determines the trajectory of the rear wheel once the initial position of the bicycle is specified. Given the initial front wheel position, the possible initial rear wheel positions constitute a circle (of radius $\l$). Given a front wheel track, the map $M:{\mathbb S}  ^1\to {\mathbb S}  ^1$ that assigns the terminal position of the bicycle to the initial one is called the {\it bicycle monodromy}. If we need to emphasize the dependence of the monodromy of the front track trajectory $F$, we write $M_F$.
In this paper, the path followed by the front wheel will   usually be closed, so $M$ is a self-map of a circle. The monodromy along a closed path depends on the starting point; another choice of this point results in a conjugated monodromy. If the  front track trajectory $F$ is not closed, we identify the initial and the terminal circles by a parallel translation; in this way, we think of the monodromy as a circle map for a non-closed path as well.

\begin{figure}[hbtp]
\centering
\includegraphics[width=2in]{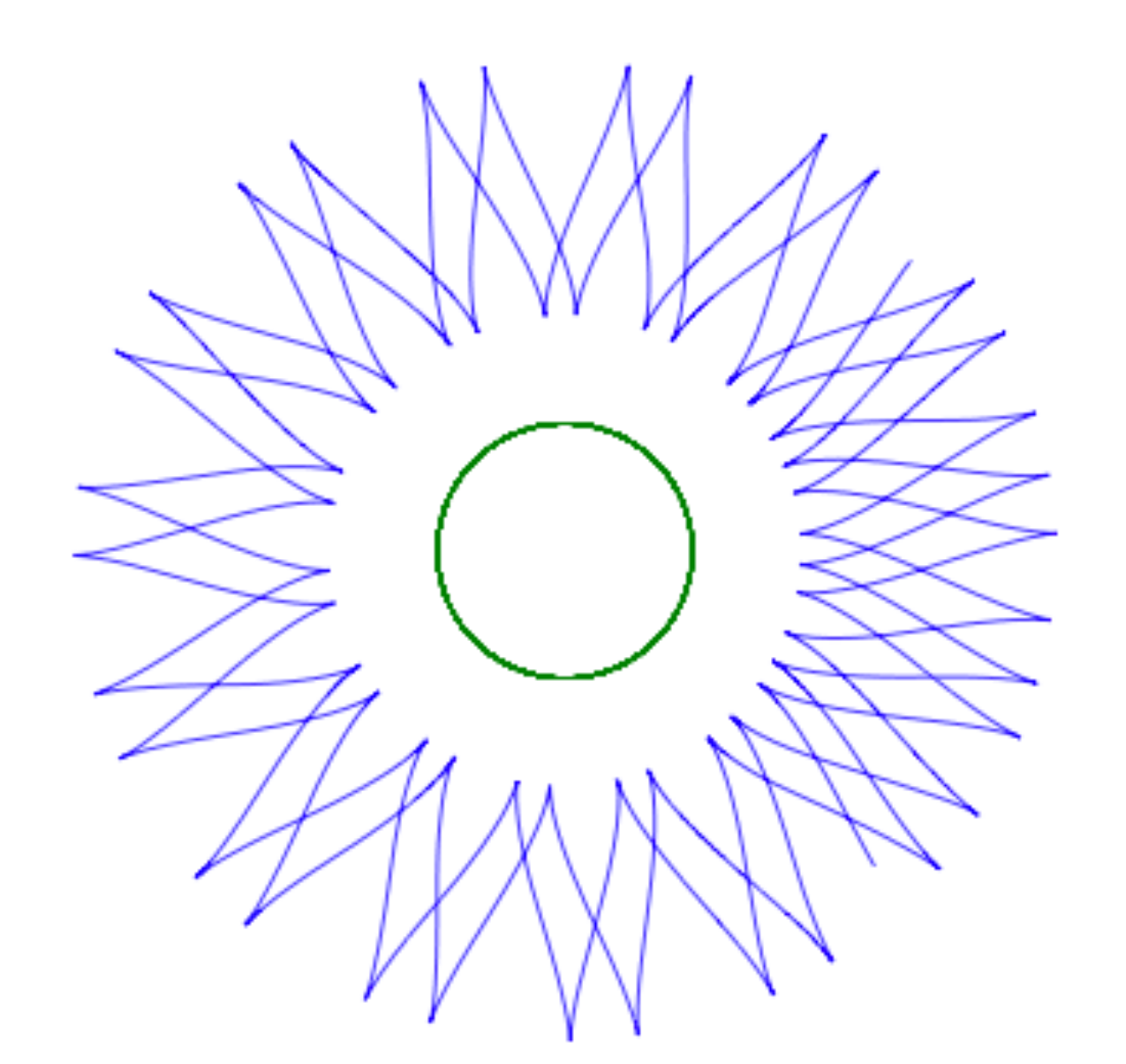}\qquad
\includegraphics[width=2in]{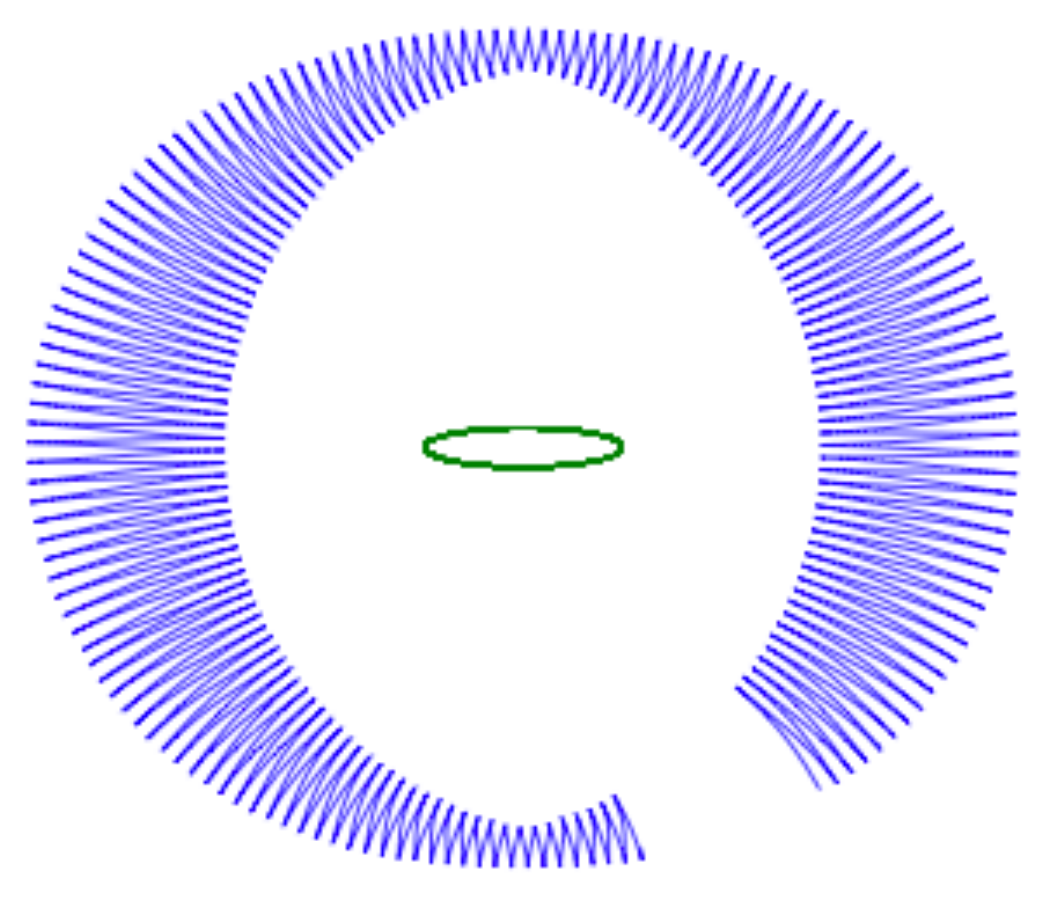}
\caption{Elliptic monodromy. The front wheel makes multiple passes around the closed curve.
}
\label{elliptic}
\end{figure}

For a closed front wheel path the circle map $M : {\mathbb S}  ^1 \to {\mathbb S}  ^1$ has two, qualitatively different, behaviors. If the length of the bicycle $\l$ is large compared to the front wheel track, one observes the behavior depicted in Figure \ref{elliptic}. In this case, the circle map $M$ is conjugate to a rotation, and it has no fixed points. This is the situation that we encountered when describing the Prytz planimeter. We refer to this behavior as elliptic.

\begin{figure}[hbtp]
\centering
\includegraphics[width=4in]{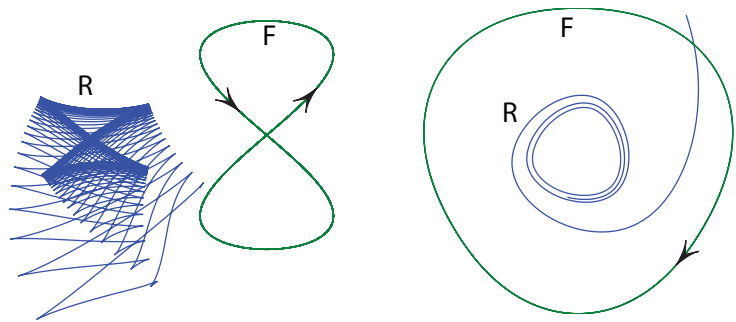}
\caption{Two examples of hyperbolic monodromy. }
\label{hyperbolic}
\end{figure}

\begin{figure}[hbtp]
\centering
\includegraphics[width=4.5in]{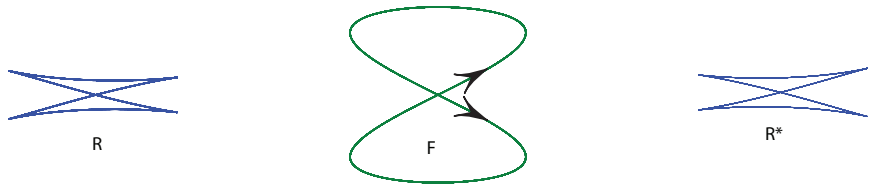}
\caption{The stable curve is on the left, and the unstable curve  on the right. If the direction of traversal of figure eight is reversed, the two curves exchange stability}
\label{eight}
\end{figure}

In contrast, when you ride a bicycle in real life, the distance the front wheel goes is generally much longer than the length of the bicycle (see Figure \ref{model}). You don't get bicycle tracks that look like those in Figure \ref{elliptic} unless you are a circus acrobat! When you make a typical round trip on a bicycle, the location of the back wheel at the end of the trip is essentially independent of its initial position: all of the possible rear-wheel trajectories are asymptotic to some particular trajectory in the family, and $M$ has an attracting fixed point, see Figure \ref{hyperbolic}. 
The map $M$ also has a repelling fixed point---it is the attracting fixed point when the bicycle runs the route in the opposite direction, see Figure \ref{eight}. This behavior is referred to as hyperbolic. For animations of this see \cite{Foote2}.

In Section~\ref{monodromy} we show that the monodromy is a M\"obius transformation of ${\mathbb S}  ^1$.
These two different behaviors are then explained by the fact that a generic fractional-linear transformation has zero or two fixed points.

A. L. Menzin was an engineer who invented a slight modification of the hatchet planimeter. More important to us, he made a conjecture in 1906 that explores the boundary between these two different types of bicycle behavior. In his own words: \cite{Men} ``If the average line across the area is long in comparison with the length of the arm, \dots\ the tractrix will approach, asymptotically, a limiting closed curve. From purely empirical observations, it seems that this effect can be obtained so long as the length of arm does not exceed the radius of a circle of area equal to the area of the base curve''. In other words,

\begin{conjecture} [Menzin] \label{Mconj}
Suppose that the path of $F$ is a simple closed curve bounding a region of area $A$. If $A > \pi\l^2$ then $M_F$ has an attracting fixed point.
\end{conjecture}

For example, if the path of $F$ is a circle of radius $r > \l$, then the chisel edge will asymptotically approach the circle of radius $\sqrt{r^2 - \l^2}$ with the same center. Similarly, if the path of $F$ is a circle of radius $r = \l$, then the chisel edge will spiral into the center. There is a qualitative difference between the two cases. If $r > \l$ then the circle of radius $\sqrt{r^2 - \l^2}$ is a periodic path attracting nearby trajectories on both sides, but if $r = \l$ then the center point is attractive only from one side and repelling from the other. As we explain below, this is the difference between hyperbolic and parabolic dynamical behavior. 

\begin{figure}[hbtp]
\centering
\includegraphics[width=3in]{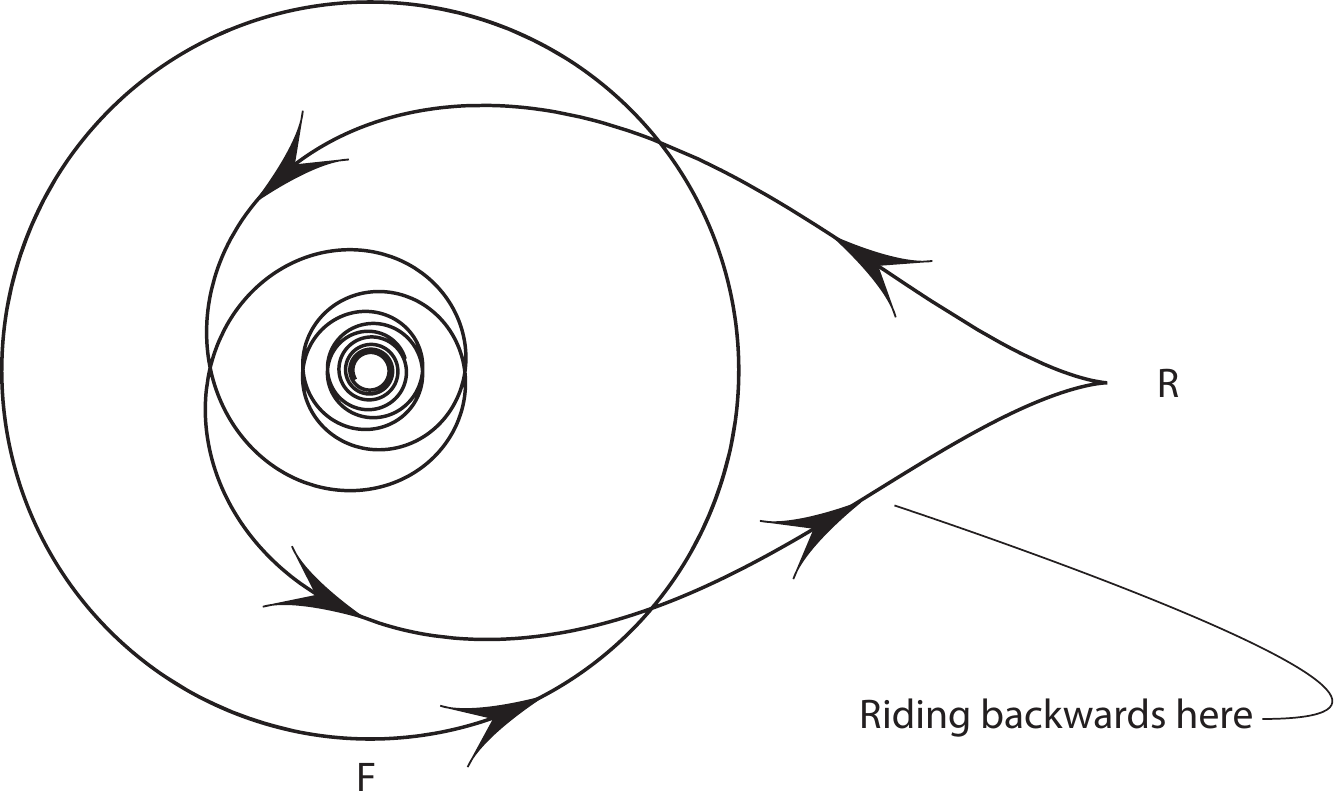}
\caption{Tractrix of a circle: $r=\l$.  
The rear wheel $R$  spirals away from the center while riding backwards; after reversing the direction at the cusp, $R$ spirals into the center. The monodromy is parabolic.}
\label{circle}
\end{figure}
 
In Section~\ref{Menzin proof} we give a proof of Menzin's Conjecture when the path of $F$ bounds a convex region.

%%%%%%%%%%%%%%%%%%%%%%%%%%%%%%%%%%%%%%%%%%%%%%%%%%%%%%%%%%%%%%%%%%%%%%%%%%%%%%%%%

\section{The Prytz planimeter and area \label{area}}

To explain how the Prytz planimeter works, let us introduce coordinates $(x,y,\th)$ on the configuration space of segments of length $\l$: the rear end (chisel edge) has coordinates $R=(x,y)$, and the direction of the segment is $\th$. The front end (tracer point) then has coordinates 
$$F=(X,Y)=(x+\l\cos\th,y+\l\sin\th).$$
The bicycle constraint is the relation $dy/dx=\tan\th$, or $\lambda=0$ where
\begin{equation} \label{1-form}
\lambda=\cos\th\ dy-\sin\th\ dx.
\end{equation}

Consider an arbitrary  motion of the segment such that its initial and  terminal positions coincide, that is, a loop in the $(x,y,\th)$-space; this motion may violate the constraint $ \lambda = 0 $.  
Denote the signed areas bounded by the closed trajectories of the rear and front ends by $A_R$ and $A_F$. We are interested in the difference $A_F-A_R$. One has:
$$
A_R=\frac{1}{2} \int xdy-ydx,\quad A_F=\frac{1}{2} \int XdY-YdX.
$$
A  computation yields:
$$
XdY-YdX = (xdy-ydx) + 2\l \lambda - \l\ d(y\cos\th - x\sin\th) +\l^2d\th.
$$
An exact differential integrates to zero over a closed curve, hence 
\begin{equation} \label{diff}
A_F-A_R = \l \int \lambda + \frac{\l^2}{2} \int d\th.
\end{equation}
The integral $\int \lambda$ measures the net violation of the bicycle constraint. In particular, $\int \lambda$ is the net signed distance that the point $R$ moves in the direction orthogonal to the segment $RF$.
The integral $\int d\th$ equals $2\pi$ times the number of turns made by the segment.

\begin{figure}[hbtp]
\centering
\includegraphics[width=2in]{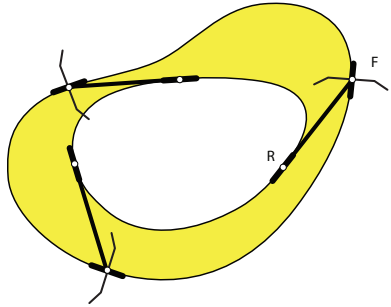}
\caption{The area between front and rear tracks is $\pi\l^2$}
\label{model}
\end{figure}

Formula (\ref{diff}) implies that the area under the tractrix (Figure \ref{tractrix}) is $\pi\l^2/2$. Indeed, the bicycle constraint is $\lambda=0$, and the moving segment turns through $180^\circ$.\footnote{Here we apply a version of (\ref{diff}) involving improper integrals, that is, integration over an infinite curve in the configuration space of segments of length $\l$.} Likewise, the area between the inner and outer tire tracks in Figure \ref{model} is $\pi\l^2$.

Now consider the hatchet planimeter in its use to measure an area. The tracer point $F$ traverses the boundary of the region being measured, starting and stopping at some base point $B$, but the chisel edge $R$ does not traverse a closed curve. Close the loop by rotating the segment  through an angle $\a$ centered at $B$, which violates the bicycle constraint. We have that $d\th$ integrates to 0 because the planimeter comes to rest in its original position without making a full rotation. Note that $\lambda = 0$ for the entire motion except the last portion, and we compute that $\lambda$ integrates to $\a\l$ which is the length of the arc followed by $R$ as it rotates about $B$ through angle $\a$. It follows from (\ref{diff}) that $A_F=\a\l^2+A_R$. Provided $A_R$ is small,  $A_F \approx \a\l^2$ is a reasonable approximation.

How good is the approximation? The error is $A_R$,  the signed area bounded by the zig-zag path of the chisel and the circular arc, Figure~\ref{planimeter}. Some starting positions are better than others, as they result in different values of $A_R$.  A suggestion made by most of the authors  is to start and stop the tracing at the centroid of the region. Of course, locating the centroid is at least as complicated as computing the area---in practice one simply makes a reasonable guess. Only Prytz \cite{Prytz} and Hill \cite{Hill} give enough mathematical details to make this rigorous (see \cite{Foote1} for a summary), and their analysis shows that even starting at the centroid does not eliminate the error entirely. They show that, for an arbitrary starting point, the error is $\origin(1/\l)$ and that, starting at the centroid, 
\begin{equation} \label{eq:smallarea}
\a\l^2 = A_F\left(1 + \frac{R^2}{2\l^2}\right) + \origin((d/\l)^3),
\end{equation} 
where $R^2$ is the mean-square distance of points in the region from the centroid and $d$ is the diameter of the region. It seems reasonable to conjecture that the full right-hand side of this formula is a weighted sum of all of the even moments of the region.

The history of the Prytz planimeter is one of humor and controversy. Other inventors, misunderstanding the mathematical nature of the error of the device, strove to improve the planimeter by adding scales or wheels that would more accurately measure the arc length $\alpha\l$ instead of the straight-line distance between the points, defeating the simple, economical design (for example, Goodman's design in Figure \ref{Prytz}). Prytz scoffed at them, writing \cite{Pr2} ``rather than use the `improved [hatchet] planimeters,' let a country blacksmith make them a copy of the original instrument.''  For more on this amusing history, see \cite{Ped, Foote1} and their references.

%%%%%%%%%%%%%%%%%%%%%%%%%%%%%%%%%%%%%%%%%%%%%%%%%%%%%%%%%%%%%

\section{Bicycle monodromy\label{monodromy}}

{\bf The monodromy is a M\"obius transformation}. The  M\"obius group on $\R $ is the group of orientation preserving isometries of the hyperbolic plane. A realization of the M\"obius group depends on the model of hyperbolic geometry. The upper half plane model identifies the hyperbolic plane with the upper half plane. The $x$-axis, complemented with a point at infinity, is the projective line $   \RP  ^1$, the absolute, or the ``circle at infinity". The group of orientation preserving isometries in this model consists of fractional-linear transformations
$$
x \mapsto \frac{ax+b}{cx+d},\quad a,b,c,d\in \R  ,\quad ad-bc>0.
$$
 
\begin{figure}[hbtp]
\centering
\includegraphics[width=2in]{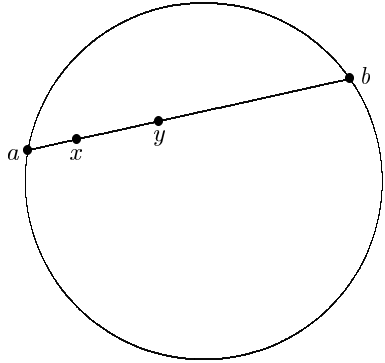}
\caption{Distance in the projective model: $d(x,y)=\frac{1}{2}\ln \frac{(a-y)(x-b)}{(a-x)(y-b)}$}
\label{cross}
\end{figure}

Another model of hyperbolic geometry is the projective (Beltrami-Cayley-Klein) model. The hyperbolic plane is represented by the interior of a disk, the lines are the chords of the disk, the distance is given by the logarithm of cross-ratio, see Figure \ref{cross}, and the the group of  isometries consists of the projective transformations of the plane that preserve the  disk. In particular, M\"obius transformations act on the boundary circle ${\mathbb S}    ^1$ of the disk. We refer to \cite{Be} for information about hyperbolic geometry.

A stereographic projection identifies the circle ${\mathbb S}  ^1$ with the projective line $\RP^1$ and conjugates the two actions, by fractional-linear transformations on $\RP^1$, and by projective transformations of the plane on ${\mathbb S}  ^1$. If $\a$ is the angular coordinate on the unit circle and $x\in \R\cup\infty$ is the coordinate on the projective line then the stereographic projection from point $(-1,0)$ is given by the formula
$x=\tan (\alpha/2)$.  

For example, the following is a 1-parameter group of projective transformations preserving the disk of radius $\l$ centered at the origin:
$$
f_t: (x,y)\mapsto \frac{\l}{\l \cosh t + x \sinh t}\ (x \cosh t + \l \sinh t, y).
$$
Since $\cosh t = 1 + \origin(t^2)$ and $\sinh t = t + \origin(t^2)$, the infinitesimal generator of this group is the vector field
\begin{equation} \label{vect}
\v w(x,y) = \frac{1}{\l}\ (\l^2-x^2, - xy)
\end{equation}
(the reader is invited to make the computations needed to verify the statements in this paragraph).

The next theorem was proved in \cite{Foote1} and extended to arbitrary dimensions in \cite{Levi}. The bicycle monodromy is a self-map of a circle of radius $\l$ which we identify with the circle at infinity in the projective model of hyperbolic geometry.

\begin{theorem} \label{Moe}
For any front wheel trajectory, the bicycle monodromy is a M\"obius transformation.
\end{theorem}

\begin{figure}[hbtp]
\centering
\includegraphics[width=1.7in]{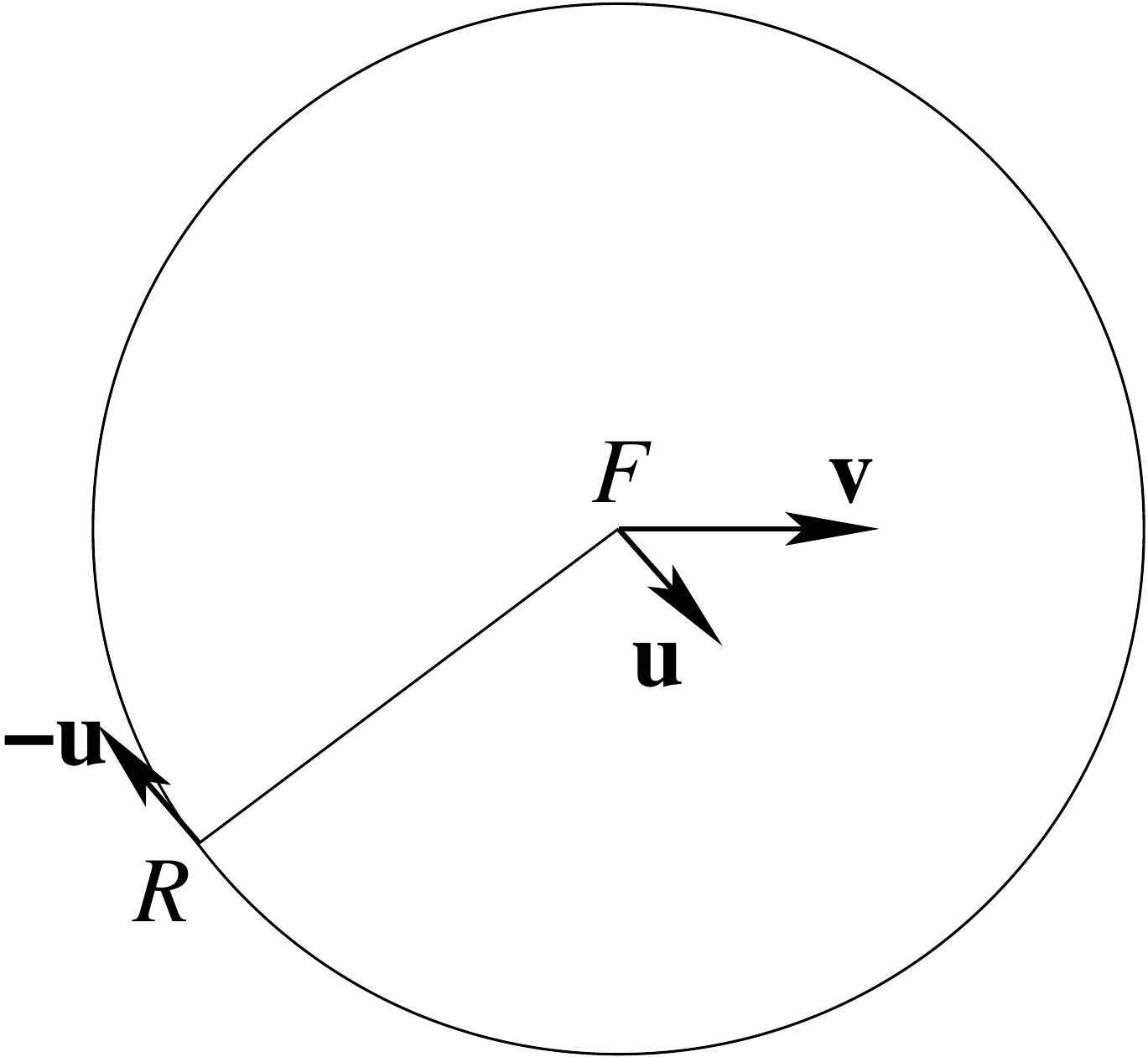}
\quad \quad \quad
\includegraphics[width=1.6in]{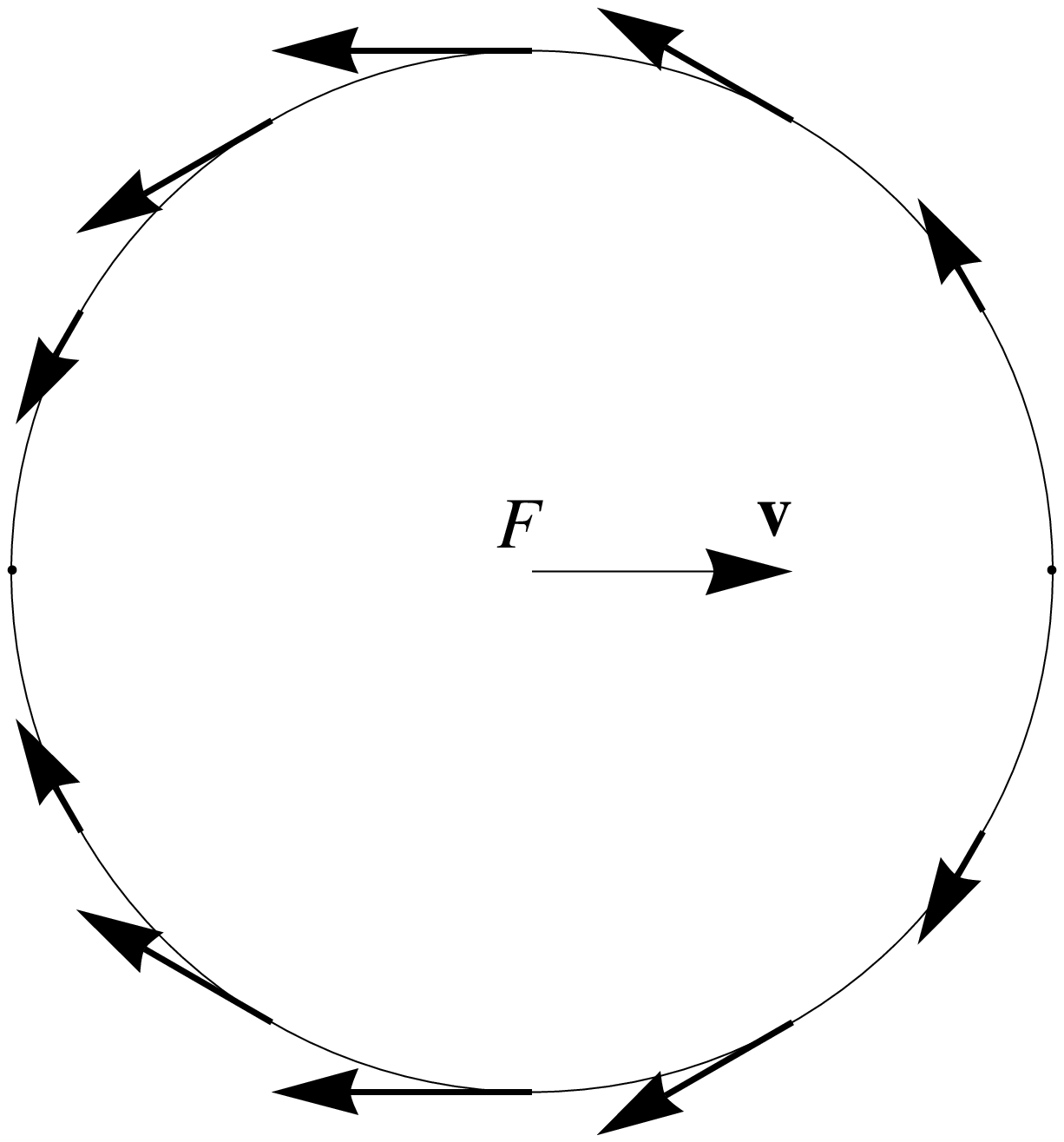}
\caption{Proof of Theorem \ref{Moe}}
\label{rel}
\end{figure}

\proof  Consider Figure \ref{rel}. We want to determine the velocity of the rear wheel $R$ relative to the front wheel $F$. Let $\v v$ be the velocity vector of $F$. Decompose this  vector into two components: the one aligned with the segment $RF$, and the perpendicular component $\v u$. If point $F$ moves along the segment $RF$ then the relative position of $R$ and $F$ does not change. On the other hand, if $F$ moves with velocity $\v u$, then $R$ moves, relative to $F$, with velocity $-\v u$. 

Thus, given a vector $\v v$, the relative velocity of every point of the circle is the negative of the projection of $\v v$ on the tangent line to the circle at this point. We have described a vector field on the circle of radius $\l$ which is an infinitesimal generator of the bicycle monodromy.

To see that this is an infinitesimal M\"obius transformation, assume (without loss of generality) that $\v v=(1,0)$. If $\a$ is the angular coordinate of point $R$ on the circle, then the negative of the projection of $\v v$ on the tangent line to the circle at point $R$ is the vector 
$\sin\a\ (\sin\a,-\cos\a)$. It remains to notice that, for $x=\l\cos\a, y=\l\sin\a$, formula (\ref{vect}) yields the same vector, scaled by $\l^2$.
\proofend

Readers familiar with differential geometry may be interested in the following interpretation: the motion of a bicycle defines a parallel translation and connection on the circle bundle over $\R^2$, and Theorem \ref{Moe} shows that the group for the connection is the M\"obius group. For details see \cite{Foote1}.

\begin{figure}[hbtp]
\centering
\includegraphics[width=1.8in]{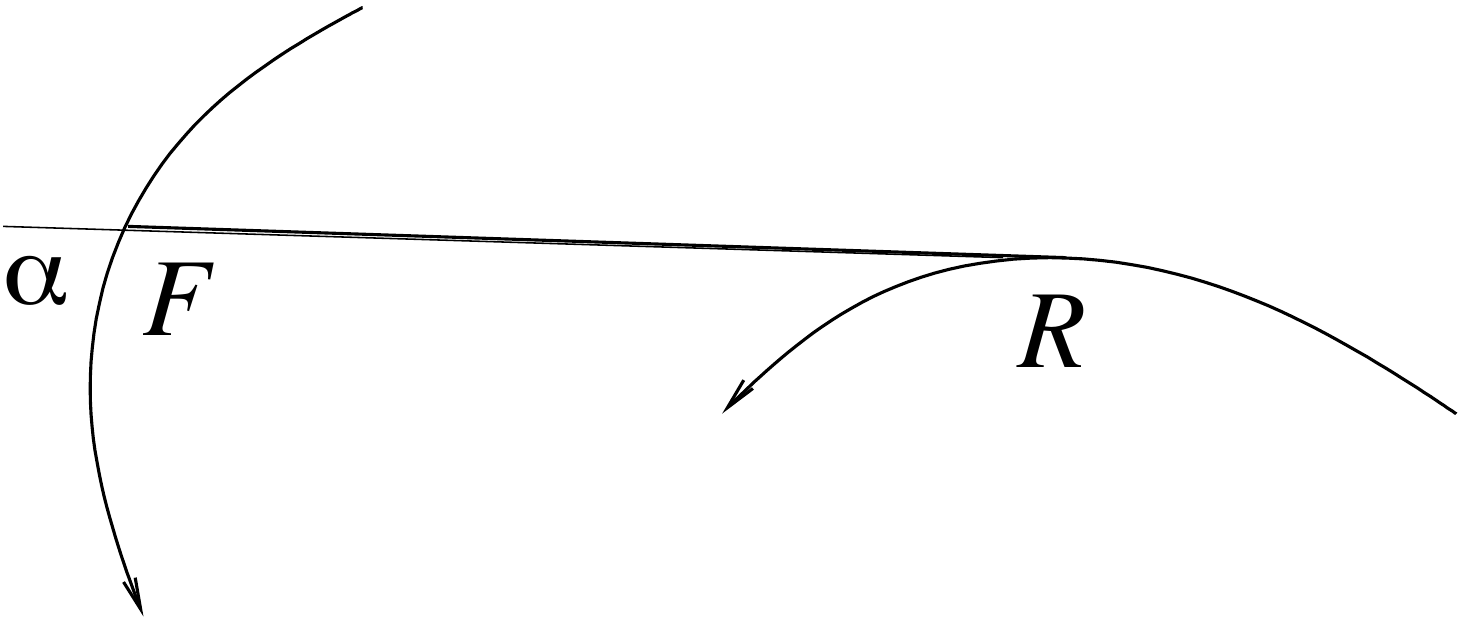}
\caption{Notation for Theorem \ref{diffeqmast}}
\label{diffeq}
\end{figure}

{\bf A differential equation}. Given a trajectory of the front wheel of the bicycle, 
%comment.6.26 its
the bike's 
 position is determined by the steering angle $\a$, see Figure \ref{diffeq}. Let $t$ be the arc length parameter along the curve $F$. The function $\a(t)$ is not arbitrary:  the bicycle constraint implies a differential equation on it. Let $k(t)$ be the curvature of the front wheel trajectory. 

Some versions of the next result appeared in \cite{AM,Dun,Finn2,Tab,Levi}, and in a different context, in \cite{CaI}.

\begin{theorem} \label{diffeqmast}
One has: 
\begin{equation} \label{master}
\frac{d\a}{dt} =k-\frac{\sin\a}{\l}.
\end{equation}
\end{theorem}

%\proof Let $T=F'$ be the unit tangent vector to the curve $F(t)$, let $S$ be the unit vector along the segment $RF$, and let $O_\a$ denote the rotation through angle $\a$, a linear transformation of  $\R^2$. Note the useful identities:
%\begin{equation} \label{ids}
%S=O_{-\a} (T),\quad (O_\a)'=\a '\ O_{\frac{\pi}{2}+\a},\quad T'=k\ O_{\frac{\pi}{2}} (T)
%\end{equation}
%(the last one holds by definition of curvature).
%
%Then one has:
%$$
%R=F-\l S=F+\l O_{\pi-\a} (T),
%$$
%and hence 
%\begin{equation} \label{R'}
%R'=T-\l \a' O_{\frac{3\pi}{2}-\a} (T) + \l k O_{\pi-\a} O_{\frac{\pi}{2}} (T) = T+\l(k-\a') O_{\frac{3\pi}{2}-\a} (T).
%\end{equation}
%The bicycle constraint implies that this vector is collinear with the vector $S$. Let $T=(t_1,t_2)$. The collinearity condition implies the equation
%$$
%\det
%\left|\begin{array}{cc}
%t_1\cos\a -t_2\sin\a+\l(k-\a')t_2,&t_1\\
%t_1\sin\a +t_2\cos\a -\l(k-\a')t_1,&t_2
%\end{array}\right| 
%=0,
%$$
%or $(t_1^2+t_2^2) (\l(k-\a')-\sin\a ) =0$. 
%Since $T\neq 0$, we arrive at  (\ref{master}).
%\proofend

\proof
Differentiating 
$
	\alpha = \arg F ^\prime -\arg RF     	
$
 by $t$ and using the definition of curvature, we get
$
	 \alpha' =    k - \omega _{RF},  
$
where $  \omega _{RF}= \frac{d}{dt}\arg RF $ is the angular velocity of $ RF $. The angular velocity is the same in all frames that do not rotate relative to each other. 
In the reference frame attached to $R$ and undergoing parallel transport, $F$ moves in a circle 
of radius $ \l $,  with speed $ v=   \sin \alpha $. Thus   $\omega _{RF} =  v/ \l = {(\sin \alpha) }/{\l} $.  
 \proofend

\begin{corollary} \label{length}
One has:
$
 R ^\prime  =  \frac{\overrightarrow{RF}}{\l}   \cos\alpha.
$
\end{corollary}

\proof Since the segment $RF$ has constant length, the speed of $R$ is the projection of the velocity of $F$ on the line $RF$, 
i.e., $\cos \alpha  $. 
%comment.6.26
And the velocity of $R$ aligns with $ \overrightarrow{RF} $. 
 \proofend

\begin{remark}
{\rm In the coordinate 
%comment.6.26 - permuted a little here. 
$x=\tan (\alpha/2)$
on the projective line, equation (\ref{master}) becomes a Riccati equation
$$
\frac{dx(t)}{dt}= \frac{1}{2} k(t) (x(t)^2+1)- \frac{1}{\l} x(t).
$$
}
\end{remark}

{\bf Signed length of the rear wheel track}. Real fractional-linear transformations come in three types: elliptic, with no fixed points; parabolic, with one neutral fixed point;  and hyperbolic, with two fixed points.
For a hyperbolic transformation, one fixed point is attractive, the other repelling, and the derivatives of the monodromy at the fixed points are  reciprocal to each other. For a parabolic transformation, the fixed point is attractive on one side and repelling on the other, and the eigenvalue at the fixed point is one. 

As we mentioned above, the rear track trajectory $R$ may have cusps. The signed length of this curve is the alternating sum of the length of its smooth pieces: the sign changes as one traverses a cusp. This signed length is the net roll of the rear wheel of the bicycle.
Consider the case when the rear track is a closed curve. Then the bicycle monodromy $M$ of the respective front track is hyperbolic or parabolic. Denote by $L$ the signed length of the rear track.

\begin{theorem} \label{eigen}
The derivatives of the monodromy at its two fixed points are equal to
$e^{\pm L/\l}.$
\end{theorem} 

\proof Let $\a(t)$ be a $T$-periodic solution of the differential equation (\ref{master}) corresponding to a fixed point of the monodromy. To find the derivative of the monodromy, consider a perturbation $\a(t)+\eps \b(t)$. Substituting into (\ref{master}) and taking the terms linear in $\eps$  yields the linearization
$$
\b'(t)=\frac{\cos\a(t)}{\l}\ \b(t).
$$
This linear equation can be solved:
$$
\b(T)=\b(0)\ e^{\int \cos\a(t)\ dt/\l}.
$$
It follows from Corollary \ref{length} that $L=\int_0^T \cos\a(t)\ dt$, hence the result.
\proofend

\begin{corollary} \label{parab}
The monodromy is parabolic if and only if $L=0$.
\end{corollary}

\proof A hyperbolic M\"obius transformation becomes parabolic when its two fixed points merge together and the derivative at this fixed point equals one.
\proofend

{\bf Interpretation of (\ref{master}): stargazing}. The differential equation (\ref{master}) describing the bicycle motion has a  curious interpretation in terms of hyperbolic geometry. Develop the front track trajectory $F(t)$ in the hyperbolic plane, that is, consider a curve $G(t) \subset H^2$ parameterized by the arc length whose curvature is the function $k(t)$. Fix a point $A$ at infinity (``an immobile star"), and let $\a(t)$ be the angle made by the line $A G(t)$ with the tangent vector $G'(t)$. 

\begin{proposition} \label{star}
One has: $\a'=k-\sin\a$.
\end{proposition}    

Thus the equation describing the motion of the unit length bicycle in the Euclidean plane also describes the the retrograde motion of the star due to motion along the curve $G$.

\begin{figure}[hbtp]
\centering
\includegraphics[width=2in]{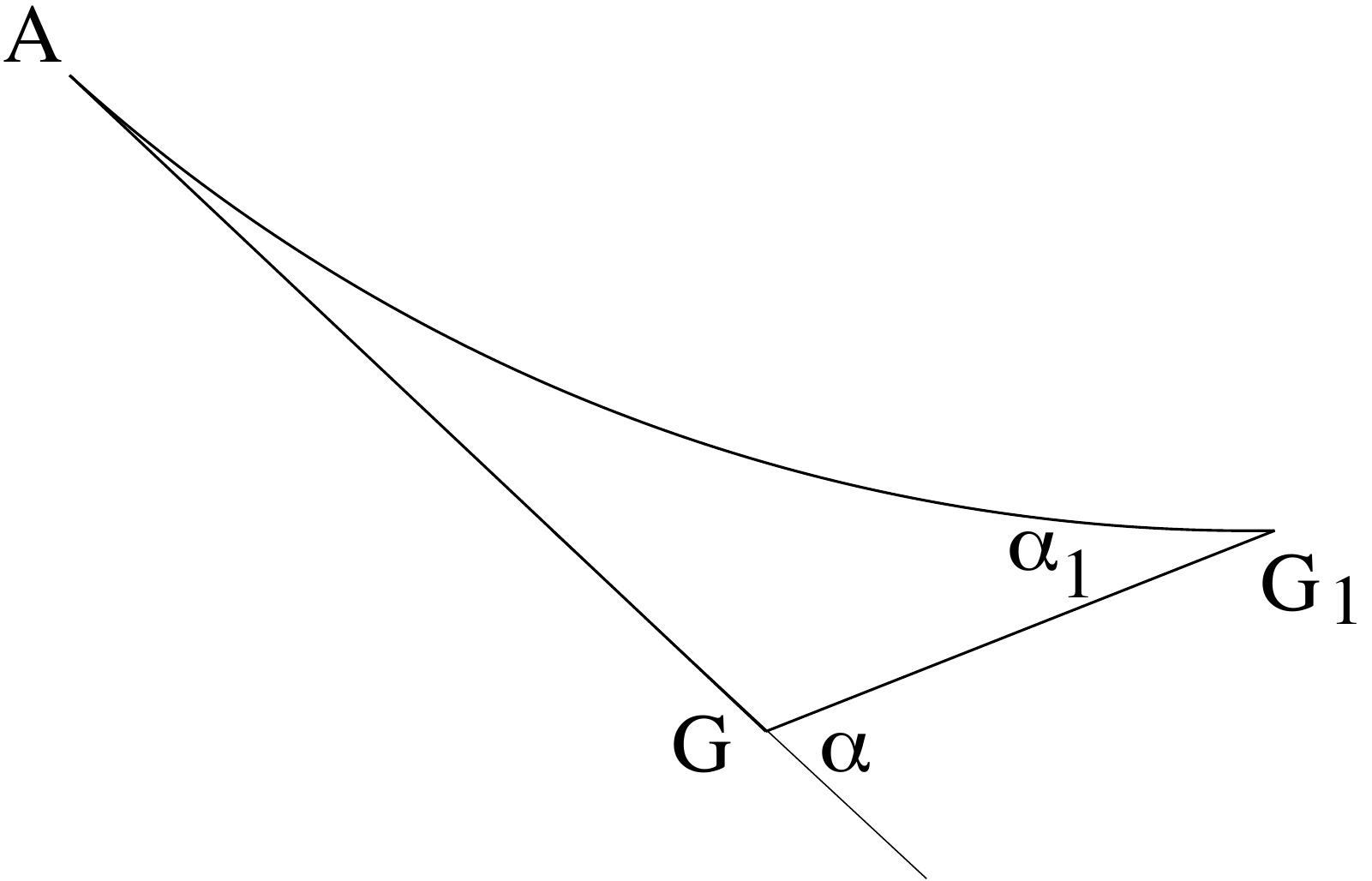}
\caption{Infinitesimal absolute triangle}
\label{tri}
\end{figure}

\proof Let $G$ and $G_1$ be infinitesimally close points on the curve. In the infinitesimally small absolute triangle $AGG_1$, the angle $A$ is zero and the side $|GG_1|=dt$. The angles of the triangle are $\pi-\a$ and $\a_1=\a+d\a$, see Figure \ref{tri}. The hyperbolic Cosine Rule (see \cite{Be}) yields:
$$
\cosh (dt)\ \sin \a\ \sin(\a+d\a) = 1-\cos \a\ \cos (\a+d\a).
$$
Expanding both sides to second order and simplifying, we get $d\alpha^2 = \sin^2\alpha\,dt^2$. Then $d\alpha = -\sin\alpha\,dt$, the sign being determined by the fact that  $d\alpha/dt < 0$ when $\sin\alpha > 0$. 

In addition, the direction of the curve changes from point $G$ to point $G_1$ by $k(t) dt$,  adding this quantity to $d\a$. Therefore $\a'=k-\sin\a$.
\proofend

We obtain a criterion for the bicycle monodromy to be the identity. Call a differentiable curve $C^1$-closed if its end points coincide and the oriented tangent lines at the end points coincide as well.

\begin{corollary} \label{ident}
The development $G\subset H^2$ of the $C^1$-closed front track $F$ is $C^1$-closed if and only if the monodromy $M_F$ for the unit length bicycle is the identity.
\end{corollary}

\proof If $G$ is $C^1$-closed then the stargazing angle $\alpha(t)\mod2\pi$ is a periodic function for all points $A$.

Conversely, assume that $G$ is not $C^1$-closed. Let $A_0$ and $A_1$ be its end points, and let $L_0$ and $L_1$ be the oriented tangent lines to $G$ at these points. If the monodromy is the identity then, for each point at infinity $X$, the lines $X A_0$ and $X A_1$ make equal angles with the curve $G$.

Denote the backward and forward intersection points of $L_0$ with the circle at infinity by $B$ and $C$. 
Then the angles  made by the lines $B A_0$ and  $C A_0$ with $G$ are zero and $\pi$ respectively, hence the lines $B A_1$ and $CA_1$ also make the angles of zero and $\pi$ with $G$, and therefore $L_1=L_0$. It follows that $A_1$ lies on $L_0$ and $A_1\neq A_0$. 

Now let $D$ be the point at infinity such that $DA_0$ is perpendicular to $L_0$. Then  $DA_1$ is not perpendicular to  $L_1\ (=L_0)$, contradicting the assumption that the monodromy is the identity. 
\proofend

For example, let $F$ be a circle of curvature $k$, traversed $p$ times, and assume that $G$ is a circle of the same curvature, traversed $q$ times. The perimeter length of $F$ is $2\pi p/k$, and that of $G$ is $2\pi q/\sqrt{k^2-1}$ (see  \cite{Be} for formulas of hyperbolic geometry). We obtain the equation
$$
\frac{p}{k}=\frac{q}{\sqrt{k^2-1}},
$$
hence $k=p/\sqrt{p^2-q^2}$. For example, $F$  can be a circle of radius $\sqrt{3}/2$ (for $p=2, q=1$), or radius $3/5$ (for $p=5, q=4$).

\begin{corollary} \label{conv}
If the front track $F$ is a $C^1$-closed convex curve (traversed once) then the bicycle monodromy $M_F$  is not the identity.
\end{corollary}

\proof Assume first that $\l=1$. 
Let $k$ be the curvature function of $F$, and assume that its hyperbolic development $G$ is also $C^1$-closed. Then  $\int_F k(t) dt =2\pi$ and, by the Gauss-Bonnet theorem in the hyperbolic plane, $\int_G k(t) dt = 2\pi +A$ where $A$ is the area bounded by $G$. This is a contradiction.

By scaling, the same  conclusion is valid for any bicycle length.
\proofend

%%%%%%%%%%%%%%%%%%%%%%%%%%%%%%%%%%%%%%%%%%%%%%%%%%%%%%%%%%%%%

\section{Proof of Menzin's conjecture \label{Menzin proof}}

In this section we prove Menzin's Conjecture \ref{Mconj} in the case when the front wheel track is convex. The property of the monodromy to be hyperbolic (and elliptic, for that matter) is ``open": a small perturbation of the curve does not affect it. Thus, without loss of generality, we assume that $F$ is a smooth strictly convex curve bounding area $A$.

The plan of the proof is to vary the length of the bicycle, from very small to very large.  We  show that
\begin{enumerate}
\item When $\l$ is small, the monodromy is hyperbolic.
\item When $\l$ is large, the monodromy is elliptic.
\item Let $\l_0$ be the smallest  length for which the monodromy becomes parabolic. Then $A\leq \pi\l_0^2$.
\end{enumerate}
Since $A>\pi\l^2$ by assumption, it follows that $\l<\l_0$, which is the hyperbolic zone by definition of $\l_0$.

\begin{figure}[hbtp]
\centering
\includegraphics[width=2.3in]{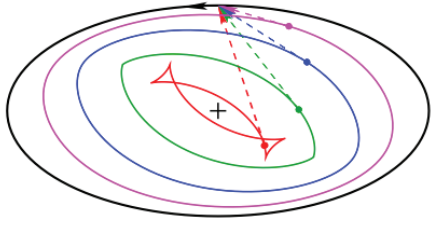}
\caption{Varying the length of the bicycle (image borrowed from \cite{Ku} with permission)}
\label{oval}
\end{figure}

See Figure \ref{oval} for a family of closed rear wheel tracks $R$ corresponding to a fixed front wheel track (outer ellipse with arrow) and various bicycle lengths. Note the change of topology of the curve $R$ as the length varies.

Let us now proceed to claims 1)--3).

  1) If $\l$ is very small, the monodromy is hyperbolic: this agrees with our experience of riding a bike whose wheel base is much smaller than the length of the path. Here is a more precise sufficient condition for hyperbolicity. 
  
 \begin{lemma} \label{suff}
 If $ \l < r$, where $r$ is the radius of the smallest osculating circle to $F$, then the monodromy  $M_F$ is hyperbolic. 
\end{lemma}   
  
  \begin{figure}[hbtp]
\centering
\includegraphics[width=2.8in]{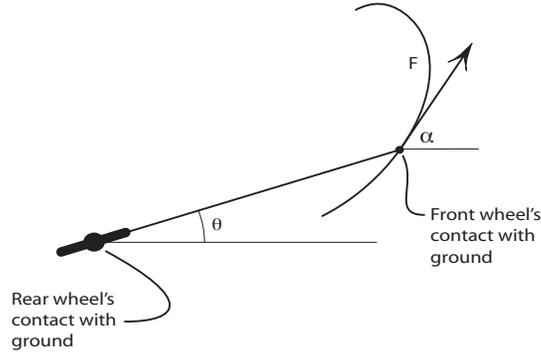}
\caption{Notations to proof of Lemma \ref{suff}.}
\label{fig:bikeode}
\end{figure}
 
\proof
 In the notations  of Figure~\ref{fig:bikeode}: 
\begin{equation} \label{eq:thetaode} 
	 \theta' = \l ^{-1} \sin(\alpha  - \theta ); 
\end{equation} 
 here $ \alpha  (t) $ is smooth   with $ \alpha  ' \geq 0 $ by the  convexity of $F$, since $ k(t) = \alpha ' (t) $ is the curvature of $F$. Let $ k_{\max}  = \max  \alpha ' (t) $. Then
\begin{equation} \label{eq:lbound}
	 \l<r= (k_{\max})^{-1}.  
\end{equation}   

\begin{figure}[hbtp]
\centering
\includegraphics[width=1.6in]{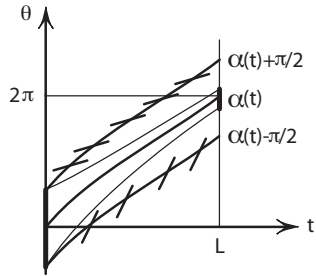}
\caption{Hyperbolicity for small $ \l $. }
\label{fig:longcurve}
\end{figure}

 Consider the strip in the $ (t, \theta ) $--plane, Figure~\ref{fig:longcurve}, given by 
 \[
	\alpha (t)  - \pi/2  \leq \theta < \alpha (t) + \pi/2 . 
\]  
This strip traps the trajectories of   (\ref{eq:thetaode}); indeed, on the lower boundary 
$ \theta = \alpha  - \pi/2 $ we have 
\[
	 \theta'(t)  = l ^{-1} \sin (\alpha  - ( \alpha  - \pi /2) )=l ^{-1}\buildrel{  (\ref{eq:lbound})}\over{>}  k_{max} \geq \frac{d}{dt} (\alpha  (t)- \pi /2).
\] 
A similar condition holds on the upper boundary of the strip, and we conclude that the segment of initial conditions 
$ [\alpha (0)- \pi/2, \alpha (0)- \pi/2] $ maps strictly into its $ 2 \pi $--translate   at $ t = L $ (the length of $F$). This implies that the monodromy map $ \theta (0) \mapsto \theta (L) $, as a map of the circle, has two fixed points -- one inside the arc from $ (\alpha (0)- \pi /2, \alpha (0)+ \pi /2) $ and another inside a complementary arc. This proves hyperbolicity. 
\proofend

2). That monodromy is elliptic for sufficiently large $\l$ was discussed in Section \ref{area}: this is what makes hatchet planimeter work. To be precise, if $F$ is convex, then $ A_F>0 $, and  we conclude, according to     (\ref{eq:smallarea}), that  if $\l$ is large enough, then the turning angle   $ 0< \alpha < 2 \pi $ for all starting points on $F$ and for all starting angles $ \theta (0) $. This proves ellipticity.

3). Now we come to the main part of the argument, the inequality $A\leq \pi\l_0^2$. Let $R_0$ be the closed rear track corresponding to the bicycle length $\l_0$. We claim that\\
(i) $R_0$ has the total rotation of $2\pi$;\\
(ii) $R_0$ is locally convex, that is, has no inflection points;\\
(iii) $A_0 \leq 0$, where $A_0$ is the area $A_0$ bounded by the curve $R_0$  given by  $(1/2) \int xdy-ydx$, as in Section \ref{area}.

Item (i) needs an explanation: even if the curve has cusps, its tangent line is well defined and continuous at all points, and we mean the rotation of this tangent line as one traverses the curve. To prove (i), notice that the total rotation depends continuously on $\l$, and that it is a multiple of $2\pi$ for a closed curve. As shown in Step 1), for small $\l$, there is a closed rear wheel curve with total rotation $2\pi$. Thus the total rotation is $2\pi$ for $R_0$ as well.

\begin{figure}[hbtp]
\centering
\includegraphics[width=3.5in]{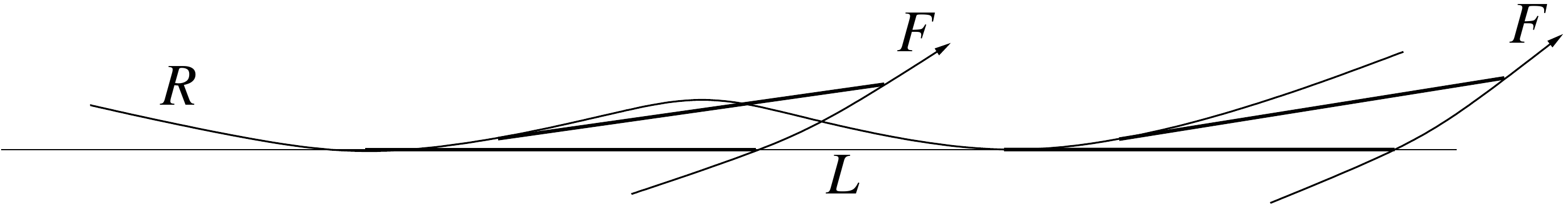}
\caption{Developing a dimple}
\label{dimple}
\end{figure}

To prove (ii), assume the opposite. Since the rear track is convex for small $\l$, non-convexity appears, as $\l$ increases, when the curvature of $R$ vanishes and then becomes negative, so that the curve develops a concavity, a ``dimple" , as shown in Figure \ref{dimple}. This yields a double tangent line $L$ to the curve $R$ which we  orient consistently with the orientation of the rear track $R$. Then the 
respective front track $F$ intersects $L$ twice, with the same intersection index (from right  to left side, in Figure \ref{dimple}).
Thus $F$ is not convex, contradicting our assumption. 

\begin{figure}[hbtp]
\centering
\includegraphics[width=2in]{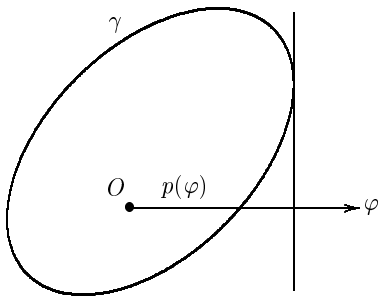}
\caption{Support function}
\label{support}
\end{figure}

To prove (iii), we use the notion of {\it support function}. Given a smooth strictly convex closed curve $\gamma$, its support function $p(\varphi)$ is the (signed) distance from the origin to the tangent line to $\gamma$, perpendicular to the direction $\varphi$, see Figure \ref{support}. The support function determines a 1-parameter family of lines, and the curve $\gamma$ is recovered as the envelope of this family. The perimeter length of the curve and the area bounded by it are given by the formulas:
\begin{equation} \label{supp}
L(\gamma)=\int_0^{2\pi} p(\varphi)\ d\varphi,\quad A(\gamma)=\frac{1}{2}\int_0^{2\pi} \left(p^2(\varphi) - p'^2(\varphi) \right) \ d\varphi,
\end{equation}
see, e.g., \cite{San}.

In fact, support functions can be used to characterize curves that are the envelopes of families of lines parameterized by their direction, that is, 1-parameter families of lines whose direction changes monotonically. As we proved,  $R_0$  is such a curve. Formulas (\ref{supp}) still apply when interpreted as signed length and signed area.

By Corollary \ref{parab}, $L(R_0)=0$. Thus we claim that 
$$
{\rm if}\ \int p(\varphi)\ d\varphi=0,\ \ {\rm then} \int p^2(\varphi)\ d\varphi\leq \int p'^2(\varphi)\ d\varphi.
$$
This is the famous Wirtinger inequality, proved by Fourier decomposition of the function $p(\varphi)$, see, e.g., \cite{Ho}.  This proves (iii).

It remains to relate the areas bounded by the rear and front tracks, that is, $A_0$ and $A$. Since the total rotation of the curve $R_0$ is $2\pi$, formula (\ref{diff}) implies:
$$
A=A_0+\pi\l_0^2\leq \pi\l_0^2,
$$
as needed. This completes the proof. 

We finish this section by describing an extension of the Menzin Conjecture and its proof to classical geometries of constant curvature, that is, to the elliptic and the hyperbolic planes \cite{Howe}. These results  were obtained  as an undergraduate research project in the REU and MASS programs at Penn State.

The bicycle constraint makes sense in both geometries, and the bicycle monodromy is still a M\"obius transformation. The spherical and hyperbolic analogs of equation (\ref{master}) are as follows:
\begin{equation} \label{sphyp}
\frac{d\alpha}{dt} = k-\cot \l\ \sin \alpha,\quad  
\frac{d\alpha}{dt} = k-\coth \l\ \sin \alpha.
\end{equation}
Here  $\cot \ell$ and $\coth \ell$ are the geodesic curvatures of the circles of radius $\l$ in ${\mathbb S}  ^2$ and $H^2$, and $k$ is the geodesic curvature of the front wheel track.

Note a curious particular case of (\ref{sphyp}) in spherical geometry: if $\l=\pi/2$ then $\cot\l=0$, and the bicycle is parallel transported along the front track $F$. If $F$ bounds area $2\pi$ then the monodromy is the identity (this is a consequence of the Gauss-Bonnet theorem).

Note also the case $\l=\infty$ in hyperbolic geometry: the second equation in (\ref{sphyp}) coincides with equation (\ref{master}) with $\l=1$.  
Since the rear end of an infinitely long bicycle in the hyperbolic plane does not move at all,  we recover Proposition \ref{star}.

An analog of the theorem proved in this section is the following result.

\begin{theorem} \label{HPZ}
In ${\mathbb S}  ^2$: if $F$ is a simple geodesically convex curve bounding area  greater than $2\pi(1-\cos \l)$  then the monodromy is hyperbolic;\\
in $H^2$: if $F$ is a simple horocyclically convex curve (i.e., having geodesic curvature greater than 1) bounding area greater than $2\pi(\cosh \l -1)$  then the monodromy is hyperbolic.
\end{theorem} 

 In both cases, the areas are those of the discs of radius $\l$.

\section{Wirtinger's inequality, Menzin's conjecture, and the isoperimetric inequality} \label{Wirtinger}

We close with two items relating our work and the isoperimetric inequality: first between our use of Wirtinger's inequality and the isoperimetric inequality for equidistant curves (wave fronts) of the rear track, and second between Menzin's conjecture and the isoperimetric inequality for the front track.

Wirtinger's inequality and the isoperimetric inequality are known to be closely related. We proved  that if a curve, given by a support function, has  zero signed length then the curve bounds non-positive signed area. We now show that this implies the classical isoperimetric inequality for convex curves. 

\begin{figure}[hbtp]
\centering
\includegraphics[width=1.8in]{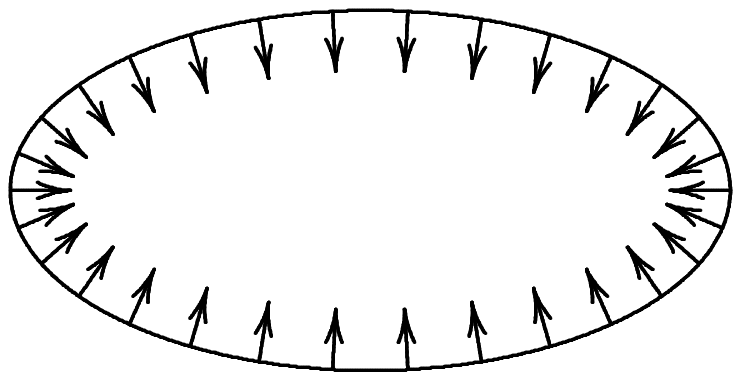}
\qquad
\includegraphics[width=2in]{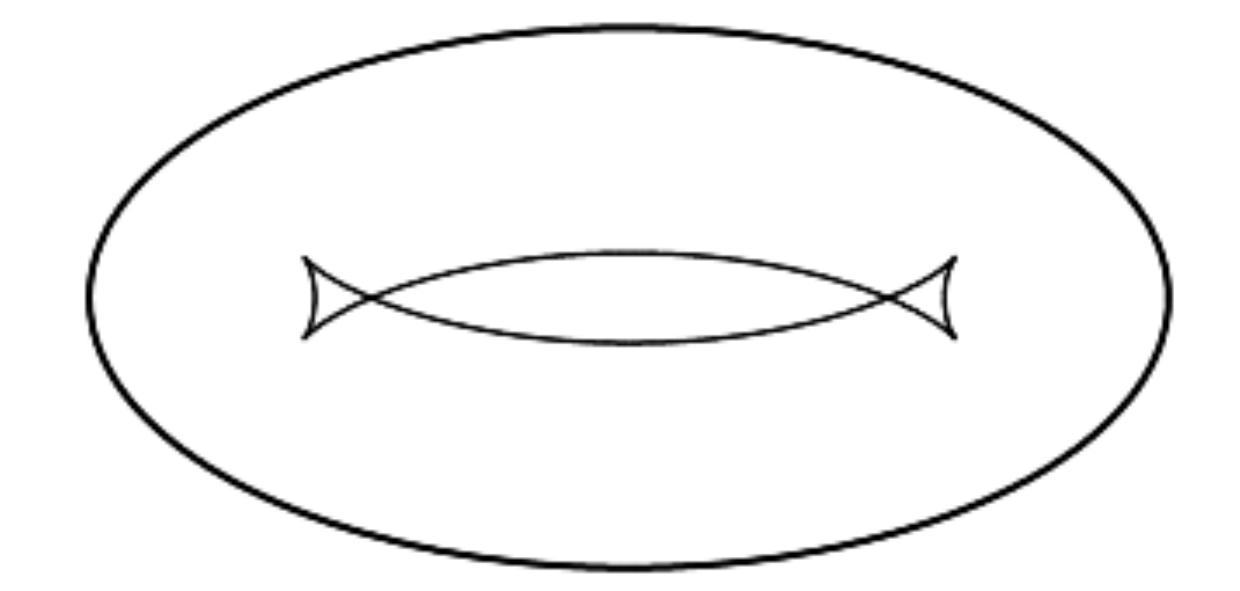}
\caption{Propagation of light inside a curve}
\label{equi}
\end{figure}

Let $\gamma$ be a closed smooth strictly convex curve, and let $L$ and $A$ be its perimeter length and the area bounded by it. Assume that $\gamma$ is a source of light, and consider propagation of light inside $\gamma$. The locus  reached by light in time $t$ is the wave front $\gamma_t$, an equidistant curve of the curve $\gamma$, see Figure \ref{equi}. The curves $\gamma_t$ are smooth and convex for small values of $t$ but later they develop cusp singularities.
The support function of $\gamma_t$ is $p_t(\varphi)=p(\varphi)-t$. It follows from formulas (\ref{supp}) that 
$$
L(\gamma_t)=L-2\pi t,\quad A(\gamma_t)=A-Lt+\pi t^2.
$$
As a consequence, the isoperimetric defect  $L^2-4\pi A$ is independent  of $t$.

Consider the moment $t$ for which $L(\gamma_t)=0$, that is, $t=L/(2\pi)$. As we proved, at this moment, 
 $A(\gamma_t)\leq 0$. Substituting $t=L/(2\pi)$ into the formula for $A(\gamma_t)$ yields the isoperimetric inequality $4\pi A\leq L^2.$ 

We conclude by another deduction of the isoperimetric inequality, for the front track, this time from Menzin's conjecture. We refer to the material of Section \ref{area}.

Suppose we bicycle around a closed path $F$ in the positive sense.
According to Menzin's conjecture, if $A_F\ge \pi\l^2$ then there is a closed rear path $R$, as in Figure \ref{model}. By equation (\ref{diff}), $A_F-A_R=\pi\l^2$, where $A_R$ is the signed area enclosed by $R$. We choose $\l$ so that $A_F = \pi\l^2$, hence $A_R=0$ (this choice of $\l$ may render the monodromy parabolic).

Let 
$$
\lambda_F = \cos\theta\,dY - \sin\theta\,dX,
$$
where $F=(X,Y)$ and $\theta$ are as in Section \ref{area}. Note that $\lambda_F = u\cdot dF$, where $u = (-\sin\theta, \cos\theta)$, and so $\lambda_F$ records the component of the motion of $F$ perpendicular to the frame of the bicycle. In particular, $\lambda_F \le dt$ where $dt$ is the length element along the front track, hence $\Lambda_F := \int \lambda_F \leq L_F$, the length of the front track.

On the other hand, it is easily shown that $\lambda_F = \lambda + \l\,d\theta$, where $\lambda$ is defined in (\ref{1-form}). The bicycle constraint is $\lambda = 0$, and so $\Lambda_F  = 2\pi\l$. Combining this with $A_F = \pi\l^2$, we have $\Lambda_F^2 = 4\pi A_F$. Thus $L_F^2\geq 4\pi A_F$, the isoperimetric inequality for the curve $F$.

Furthermore, the equality implies that $\lambda_F = ds$, that is, the front wheel only moves 
in the direction perpendicular to the bicycle frame. In this case the rear wheel does not move at all, 
and the front wheel describes a circle of radius $\l$. Thus the isoperimetric inequality is an
equality only if front wheel path is a circle.

It's worth noting that the inequality
$$
L_F^2 \ge \Lambda_F^2 = 4\pi(A_F - A_R)
$$
 holds whenever there is a closed rear wheel path. In the proof of Menzin's conjecture, we identified a closed rear wheel path with signed area $A_0 \le 0$. In cases where $A_0 < 0$, we get a positive lower bound on the isoperimetric defect:
$$
L_F^2 - 4\pi A_F \ge \Lambda_F^2 - 4\pi A_F = -4\pi A_0.
$$

\bigskip
{\bf Acknowledgments}. We are grateful to many a mathematician for discussions of  ``bicycle mathematics". In particular, it is a pleasure to acknowledge interesting discussions with A. Kurnosenko, J. Langer, R. Perline. 
 S. T. was partially supported by the Simons Foundation grant No 209361 and by the NSF grant DMS-1105442; M.L. was supported by the NSF grant DMS-0605878.

\end{document}